\documentclass[12pt]{article}
\usepackage{a4wide}

\usepackage{amsmath,amssymb}
\usepackage{bbm}
\usepackage{graphicx}

\newcommand{\cv}{\mathbf{c}}
\newcommand{\uv}{\mathbf{u}}
\newcommand{\rv}{\mathbf{r}}
\newcommand{\ev}{\mathbf{e}}
\newcommand{\Ev}{\mathbf{E}}
\newcommand{\Dv}{\mathbf{D}}

\newcommand{\dv}{\mathbf{d}}
\newcommand{\vv}{\mathbf{v}}
\newcommand{\fv}{\mathbf{f}}
\newcommand{\xv}{\mathbf{x}}
\newcommand{\Xv}{\mathbf{X}}
\newcommand{\nv}{\mathbf{n}}

\newcommand{\Dtm}{\mathsf{D}}
\newcommand{\Etm}{\mathsf{E}}
\newcommand{\Ctm}{\mathsf{C}}

\newcommand{\dtm}{\mathsf{d}}
\newcommand{\etm}{\mathsf{e}}
\newcommand{\ctm}{\mathsf{c}}

\newcommand{\Ft}{\mathbf{F}}
\newcommand{\Ct}{\mathbf{C}}
\newcommand{\St}{\mathbf{S}}
\newcommand{\It}{\mathbf{I}}
\newcommand{\Tt}{\mathbf{T}}

\newcommand{\sigmat}{\boldsymbol{\sigma}}
\newcommand{\taut}{\boldsymbol{\tau}}

\newcommand{\epst}{\boldsymbol{\varepsilon}}

\newcommand{\Body}{\Omega}
\newcommand{\body}{\omega}
\newcommand{\opdiv}{\operatorname{div}}
\newcommand{\opcurl}{\operatorname{curl}}
\newcommand{\sym}{\operatorname{sym}}
\newcommand{\skw}{\operatorname{skw}}

\usepackage{color}

\begin{document}
\title{Large deformation mixed finite elements for smart structures} 
\author{Astrid S. Pechstein}

\date{\today}
\maketitle
\begin{abstract}{
Recently, "Tangential Displacement Normal Normal Stress" (TDNNS) elements were introduced for small-deformation piezoelectric structures. Benefits of these elements are that they are free from shear locking in thin structures and volume locking for nearly incompressible materials. We extend these elements to the large deformation case for electro-active polymers in the framework of an updated Lagrangian method. We observe that convergence does not deteriorate as the material becomes nearly incompressible with growing Lam\'e parameter $\lambda$, and that the discretization of slender structures by flat volume elements is feasible.
The elements are freely available in the open source software package Netgen/NGSolve.}
\end{abstract}

\section{Introduction}

Smart materials are state of the art in realizing mechatronic applications, both as actuators and sensors. Applications range from active noise or vibration control, to energy harvesting, haptic feedback, and precision systems. In the small deformation range, piezoceramics such as PZT are the method of choice in smart structure applications. More recently, intensive research on electro-active polymers (EAPs) enlarges the set of suitable smart materials. These polymers undergo large deformations under electric loading.

In the present contribution, we are concerned with the simulation of large deformation nonlinear electro-elasticity. For an introduction to nonlinear electro-elasticity we refer to the monograph by Maugin \cite{Maugin:1988}. Electro-active polymers have been characterized in the framework of nonlinear electro-elasticity or hyperelasticity in \cite{VoltairasFotiadisMassalas:2003,DorfmannOgden:2005,SuoZhaoGreene:2008}. Bustamante \emph{et al.}~\cite{BustamanteDorfmannOgden:2009}  derived different variational principles based on the different formulations. Stability of the incremental equations is addressed in \cite{DorfmannOgden:2010}, while invariants and restrictions of the constitutive laws are treated in \cite{BustamanteMerodio:2011}. A multiplicative splitting of mechanic and electric deformation is suggested by Skatulla~\emph{et al.} \cite{SkatullaSansourArockiarajan:2012}.  Different multiplicative splittings that are micromechanically motivated are described by Z\"ah and Miehe \cite{ZaehMiehe:2015}.

Based on the description by Dorfmann and Ogden \cite{DorfmannOgden:2005}, Vu \emph{et al.} \cite{VuSteinmannPossart:2007} proposed a finite element formulation for electro-active polymers. They developed volume finite elements with four degrees of freedom per node---three for the displacements and one for the electric potential. They proved that the element is well-suited for the large deformation case. 

However, it is well known that nodal finite elements suffer from several types of locking, such as volume locking for nearly incompressible materials, or bending and shear locking for flat elements. Mixed methods have been found to avoid such locking phenomena \cite{BoffiBrezziFortin:2013}. The price for stability and accuracy is the introduction of additional unknowns, e.g.\ stresses, pressure, deformation gradient, strain or dielectric displacement. We shortly review on different solution strategies proposed in the literature, where we differ between ``flat'' elements, such as plate, shell or solid shell elements, and volume elements that are not (or only weakly) sensitive to large aspect ratios.

For piezoelectric materials in the small-strain regime, we mention hybrid stress solid shell  elements by Sze~\emph{et al.}~\cite{SzeYaoYi:2000}, where stresses are added explicitely. Klinkel and Wagner \cite{KlinkelWagner:2006} proposed a solid shell element based on a Hu-Washizu formulation with six independent fields. While this element uses a linear constitutive law for piezoelasticity, a similar element employing a nonlinear constitutive law suitable for dielectic elastomers is described by Klinkel~\emph{et al.}~\cite{KlinkelZweckerMueller:2013}. Krommer \emph{et al.} \cite{KrommerVetyukovStaudigl:2016} reduce the material laws from three to two dimensions analytically, assuming plane stress state.

Small-strain volumetric elements for piezoelectric materials were designed by Sze and Pan \cite{SzePan:1999}. These elements use assumed stesses and electric displacements, which are condensed at the element level. More recently, Ortigosa and Gil \cite{OrtigosaGil:2016} introduced volume elements suitable for large deformations. These elements rely on a multi-variable convex potential and include five additional unknown fields.

In the present contribution, a finite element method is introduced where displacement, total stress and electric potential are the unknown fields. It is based on the ``Tangential Displacement Normal Normal Stress'' (TDNNS) method, which is introduced and analyzed for small, purely elastic deformations in \cite{PechsteinSchoeberl:2011,PechsteinSchoeberl:2016}. In \cite{PechsteinSchoeberl:2012} it was shown that prismatic volume elements are free from shear locking and do not suffer from large aspect ratios. Additionally, the elements are suitable for the description of nearly incompressible materials \cite{Sinwel:2009} if a consistent stabilization term is added. An extension to small-deformation applications for piezoelectric solids can be found in \cite{PechsteinMeindlhumerHumer:2018}.

In our deduction, we start from any standard energy formulation such as the Neo-Hookean-type free energy for electro-active polymers used in \cite{VuSteinmannPossart:2007}.  We derive a formulation close to the updated Lagrangian setting \cite{BatheRammWilson:1975}. The updated Lagrangian method is in general equivalent to common Lagrangian methods, but has been shown to benefit from the possibility of adaptive remeshing \cite{LegerEtal:2014,LegerPepin:2016}. A consistent variational formulation for the linearized update equation around some intermediate configuration is derived. Displacement update, total stress and electric potential update in intermediate configuration are the unknown fields of the method. 

The paper is organized as follows: In Section~\ref{sec:notation} we introduce notation for tensor calculus. Section~\ref{sec:largedef} addresses large deformation electro-elasticity, having in view the updated Lagrangian setting. In Section~\ref{sec:mixed} the finite elements as known from the linear problem are recalled briefly, and the mixed formulation is derived analytically. Finally, in Section~\ref{sec:numerics} we present numerical results.

\subsection{Notations}\label{sec:notation}

Throughout this paper, we will use a direct tensor notation, where vectors are interpreted as first-rank tensors, and the dot $\cdot$ denotes the contraction, e.g.
\begin{align}
\sigmat \cdot \taut &= \sigma_{ik} \tau_{kj}, & \sigmat : \taut &= \sigma_{ij} \tau_{ji}.
\end{align}
Above and wherever it seems advantageous, we switch to component-wise notation using Einstein's summation convention. Symmetric and skew-symmetric part of a second order tensor are defined as 
\begin{align}
\sym(\sigmat) & = \frac{1}{2} (\sigmat + \sigmat^T), & \skw(\sigmat) &= \frac{1}{2}(\sigmat - \sigmat^T).
\end{align}

For some domain $D$, let $\partial D$ denote its boundary. We denote the outward unit normal by $\nv$. The normal and tangential component of a vector field $\vv$ are defined as
$v_n = \vv \cdot \nv$ and $\vv_t = \vv - v_n \nv$, respectively.
For a  tensor field $\sigmat$ we define the normal vector $\sigmat_n = \sigmat \cdot \nv$ and use its normal and tangential components $\sigma_{nn}$ and $\sigmat_{nt}$.

Derivatives with respect to spatial coordinates will be addressed by the common nabla, divergence and curl operators. As this is important in our large-deformation setting, we will always indicate the coordinate system in which the differentiation is performed.
The nabla or gradient operator with respect to coordinates $\xv$ is defined as
\begin{align}
\nabla_\xv &= \frac{\partial}{\partial x_i}.
\end{align}
We use the divergence operator acting on a vector field $\vv$ or a tensor field $\sigmat$,
\begin{align}
\opdiv_\xv \vv &= \nabla_\xv \cdot \vv = \frac{\partial v_i}{\partial x_i}, \qquad\text{and}&
\opdiv_\xv \sigmat &= \nabla_\xv \cdot \sigmat = \frac{\partial \sigma_{ij}}{\partial x_i} .
\end{align}
The curl operator of a three-dimensional vector field is defined in the standard way.
We note that the curl operator is equivalent to the skew-symmetric part of the gradient, as this will be a principal ingredient in the derivation of the mixed method,
\begin{align}
\nabla^{skw}_\xv \vv  = \skw(\nabla_\xv \vv) &= \begin{array}{cc}
\frac12 \left[ \begin{array}{ccc} 0 & -c_3 & c_2 \\ c_3 & 0 & -c_1 \\ -c_2 & c_1 & 0 \end{array}\right] & \text{ for }\cv = \opcurl_\xv \vv.
\end{array}.
\label{curlskw}
\end{align}
A similar relation  holds also for the two-dimensional setting.

\section{Electromechanics in the large deformation setting} \label{sec:largedef}

In this section, we introduce our notation for the large deformation problem, such as material coordinates of the reference configuration and spatial coordinates of the deformed configuration. The basic electro-mechanical fields are discussed, as well as their material and balance laws in Lagrangian and spatial configuration.

\subsection{Material and spatial configuration}

Let $\Body \subset \mathbb{R}^d, d=2,3$ denote the body of interest in undeformed configuration. The motion, or deformation of $\Body$ is described by the mapping $\rv: \Body \to \body$, where $\body$ is the current or deformed configuration in space. Each material point $\Xv \in \Body$ is mapped to a spatial point $\xv \in \body$ via
\begin{equation}
\xv = \rv(\Xv) = \Xv + \uv(\Xv),
\end{equation}
with $\uv$ the displacement. The deformation gradient $\Ft$ is the material derivative of $\rv$,
\begin{align}
\Ft &= \nabla_\Xv \rv = \It + \nabla_\Xv \uv.
\end{align}
We define the right Cauchy-Green tensor $\Ct = \Ft^T \cdot \Ft$, the strain $\St = \frac12 (\Ct - \It)$ and the Jacobian $J = \det \Ft$.
We further use the symmetic second Piola-Kirchhoff stress $\Tt$, which satisfies the mechanical equilibrium condition
\begin{align}
\opdiv_\Xv ((\Ft \cdot\Tt)^T) &= 0 &&\text{in } \Body.
\end{align}
The analogous equation in spatial configuration uses the symmetric total stress $\taut$,
\begin{align}
\opdiv_\xv \taut &= 0 &&\text{in $\body$, with} & \taut &= \frac{1}{J} \Ft\cdot \Tt\cdot \Ft^T.
\end{align}

Throughout the following, we assume that the assumptions of electrostatics are applicable, which is usually the case as the speed of light exceeds the speed of sound by orders of magnitude. In this case, the electric field  $\ev$ and the dielectric displacement $\dv$ are governed by Gauss's law for electricity and Faraday's law of induction,
\begin{align}
\opdiv_\xv \dv &= 0 &\text{and}&& \opcurl_\xv \ev &= 0 & \text{in } \body.
\end{align}
Note that both laws are formulated in current, not in reference configuration. The electric field is then a gradient field, using the electric potential $\phi$,
\begin{align}
\ev &= -\nabla_\xv \phi.
\end{align}

These definition can be transformed to material configuration. As one can verify using basic calculus, the material electric field $\Ev$ is linked to the spatial electric field $\ev$ via
\begin{align}
\Ev &= \nabla_\Xv \phi = \Ft^T\cdot \ev.
\end{align}
The dielectric displacement is transformed by the Piola transformation, and satisfies
\begin{align}
 \Dv &= J \Ft^{-1}\cdot \dv, & \opdiv_\Xv \Dv &= 0.
\end{align}

Usually, the material law provides the connection between strain, electric field, stress and dielectric displacement in material configuration. Throughout the following, we assume that we are given a free energy function $\Psi(\Ct, \Ev)$. Then, stress and dielectric displacement are defined as
\begin{align}
\Tt &= 2 \frac{\partial \Psi}{\partial \Ct}, & \Dv &= - \frac{\partial \Psi}{\partial \Ev}. \label{eq:TD}
\end{align}

\subsection{Linearization at a given state}

In the following, we gather the linearized equations at some given state $(\uv_0, \phi_0)$. In our finite element methods, all computations are done in deformed (or ``updated Lagrangian'' \cite{BatheRammWilson:1975}) configuration. With $\rv_0$ the deformation map corresponding to $\uv_0$, we define $\omega_0$ as the image of $\Omega$ under $\rv_0$, i.e. $\rv_0: \Omega \to \omega_0$. We are interested in linearizing the equations for the displacement and electric potential updates $(\uv_\Delta, \phi_\Delta)$, which are defined on $\omega_0$ such that
\begin{align}
\uv(\Xv) &= \uv_0(\Xv) + \uv_\Delta(\rv_0(\Xv)),& \phi(\Xv) &= \phi_0(\Xv) + \phi_\Delta(\rv_0(\Xv)).
\end{align}
We assume the displacement update to be small, and identify $\omega_0$ with $\omega$ and $\xv$ with $\rv_0(\Xv)$.
The deformation gradient decomposes multiplicatively,
\begin{align}
\Ft &= \Ft_\Delta \cdot\Ft_0, & \text{with } \Ft_0 &= \It + \nabla_\Xv(\uv_0),\quad \Ft_\Delta = \It + \nabla_\xv (\uv_\Delta).
\end{align}
The linearized update of the right Cauchy-Green tensor is well known to be
\begin{align}
\Ct - \Ct_0 &= \Ft_0^T\cdot \Ft_\Delta^T \cdot\Ft_\Delta\cdot \Ft_0 - \Ft_0^T \cdot\Ft_0 = 2 \Ft_0^T \cdot\epst_{\xv,\Delta} \cdot\Ft_0, \label{linearCupdate}
\end{align}
with $\epst_{\Delta} = \sym(\nabla_\xv \uv_\Delta)$. The update of the material electric field directly computes as
\begin{align}
\Ev - \Ev_0 &= - \nabla_\Xv \phi + \nabla_\Xv \phi_0 = - \nabla_\Xv \phi_\Delta. \label{Eupdate}
\end{align}

When linearizing the material law \eqref{eq:TD} in Lagrangian configuration, we use \eqref{linearCupdate} and \eqref{Eupdate} to obtain the material moduli $\Ctm, \Dtm$ and $\Etm$ as certain derivatives of $\Psi$ with respect to $\Ct$ and $\Ev$, 
\begin{align}
\Tt &= 2 \frac{\partial \Psi}{\partial \Ct}|_{(\uv_0, \phi_0)} + 2 \frac{\partial^2 \Psi}{\partial \Ct^2}|_{(\uv_0, \phi_0)}: (\Ct - \Ct_0) + 2 \frac{\partial^2 \Psi}{\partial \Ct\partial \Ev}|_{(\uv_0, \phi_0)}\cdot(\Ev - \Ev_0)\label{mat1}\\
&= \overline{\Tt}_0 + \Ctm: (\Ft_0^T \cdot \epst_{\Delta} \cdot\Ft_0) - \Etm \cdot\nabla_\Xv \phi_\Delta, \nonumber\\
\Dv &= - \frac{\partial \Psi}{\partial \Ev} - \frac{\partial^2 \Psi}{\partial \Ev\partial \Ct}|_{(\uv_0, \phi_0)} :(\Ct - \Ct_0) - \frac{\partial^2 \Psi}{\partial \Ev^2}|_{(\uv_0, \phi_0)} \cdot(\Ev - \Ev_0) \label{mat3}\\
&= \overline{\Dv}_0 - \Etm: (\Ft_0^T \cdot\epst_{\Delta}\cdot \Ft_0) + \Dtm \cdot\nabla_\Xv \phi_\Delta. \nonumber
\end{align}
Above, we also implicitly defined 
\begin{align}
\overline{\Tt}_0 &= 2 \frac{\partial \Psi}{\partial \Ct}|_{(\uv_0, \phi_0)} , & \overline{\Dv}_0 &= - \frac{\partial \Psi}{\partial \Ev}.
\end{align}
By transformation we gain the spatial moduli $\ctm$, $\etm$ and $\dtm$,
\begin{align}
\begin{aligned}
\ctm_{ijkl} &= \frac{1}{J_0} \Ctm_{IJKL} F_{0,iI} F_{0,jJ} F_{0,kK} F_{0,lL},\\
\etm_{ikl} &= \frac{1}{J_0} \Etm_{IKL} F_{0,iI} F_{0,kK} F_{0,lL},&
\dtm_{ij} &= \frac{1}{J_0} \Dtm_{IJ} F_{0,iI} F_{0,jJ},
\end{aligned}
\end{align}
and the total stress and dielectric displacement in spatial configuration
\begin{align}
\overline{\taut}_0 & = \frac{1}{J_0} \Ft_0 \cdot \overline{\Tt}_0 \Ft_0^T, & \overline{\dv}_0 &= \frac{1}{J_0} \Ft_0 \cdot \overline{\Dv}_0.
\end{align}
When including geometric stiffening effects, the governing equation in deformed configuration $\omega_0$ read
\begin{align}
\Ft_\Delta^{-1}\cdot \taut &= \overline{\taut}_0 + \ctm  : \epst_{\Delta}  - \etm\cdot \nabla_\xv \phi_\Delta \label{pde1}\\
\dv &= \overline{\dv}_0 - \etm\cdot \epst_\Delta + \dtm \cdot\nabla_\xv \phi_\Delta, \label{pde2}\\
\opdiv_\xv \taut^T &= 0, \label{pde3}\\
\opdiv_\xv \dv &= 0. \label{pde4}
\end{align}
Note that the left-hand side of \eqref{pde1} still contains some nonlinearity, which will be treated in Section~\ref{sec:mixedfe}.

\section{A mixed variational formulation} \label{sec:mixed}

In \cite{PechsteinMeindlhumerHumer:2018}, we proposed mixed finite elements for linear piezoelectric materials. In the current contribution, we will use these elements for the total stress, the displacement and electric potential updates. We will shortly present the finite elements in the linear setting, before adapting the equations to the nonlinear case.

Throughout the remainder of the paper, we assume $\mathcal{T}(\Omega) = \{T\}$ to be a finite element mesh of the reference configuration. By $\mathcal{T}(\omega_0) = \{t_0\}$ we mean the finite element mesh transformed by $\rv_0$. The method works for hybrid meshes consisting of triangular and quadrilateral in 2D, or tetrahedral, hexahedral and prismatic elements in 3D.

\subsection{Mixed elements for the linear problem}

The proposed finite elements have non-standard degrees of freedom, namely the tangential component of the displacement and the normal component of the stress vector on element interfaces. This implies that the displacement field is not continuous, and that gaps in normal direction may open up between elements. The normal stresses then control these gaps. The elements were introduced and analyzed in the small-strain elastic setting in \cite{PechsteinSchoeberl:2011,PechsteinSchoeberl:2012,PechsteinSchoeberl:2016}, and applied to piezoelectric materials in \cite{PechsteinMeindlhumerHumer:2018}. 

We shortly present the mixed method in the small-strain piezoelectric case, using the (then constant) material moduli from equations \eqref{mat1} and \eqref{mat3}. For solution and virtual quantities, we use the spaces (assuming simplicial elements, otherwise the polynomial degrees vary for the different components, see \cite{Zaglmayr:2006,PechsteinSchoeberl:2012})
\begin{align}
\uv, \delta \uv &\in \mathbf{V}(\Omega) := \{ \vv: \vv|_T \in P^{k}(T), \vv_t \text{ cont.}, \vv_t = 0\text{ on }\Gamma_{fix}\},
\label{spaceX1}\\
\taut, \delta \taut & \in \boldsymbol{\Sigma}(\Omega) := \{ \sigmat: \sigmat \text{ symm}, \sigmat|_T \in P^{k}(T), \sigma_{nn} \text{ cont.}, \sigma_{nn} = 0\text{ on }\Gamma_{free}\},
\label{spaceX2}\\
\phi, \delta \phi & \in W(\Omega):= \{ w: w|_T \in P^{k+1}(T), w \text{ cont.}, w = 0 \text{ on }\Gamma_E\}.
\label{spaceX3}
\end{align}
We note that, for a displacement field $\uv \in \mathbf{V}(\Omega)$, neither displacement gradient nor linearized strain tensor exist in $L^2$ sense. Indeed, the strain is a distribution, but work pairs $\langle \taut, \uv\rangle_\Omega$ can be evaluated for $\taut \in \boldsymbol{\Sigma}(\Omega)$ by
\begin{align}
\langle \taut, \epst(\uv)\rangle_\Omega &= \sum_{T \in \mathcal{T}} \Big( \int_T \taut : \epst(\uv)\, d\Xv - \int_{\partial T} \tau_{nn} u_n\, dS_X \Big) \label{epsdiv1}\\
&= - \sum_{T \in \mathcal{T}} \Big( \int_T \opdiv_\Xv \taut \cdot \uv\, d\Xv - \int_{\partial T} \taut_{nt}\cdot \uv_t\, dS_X \Big) = -\langle \opdiv_\Xv \taut, \uv\rangle_\Omega. \label{epsdiv2}
\end{align}
For further information and a thorough analysis we refer the interested reader to \cite{PechsteinSchoeberl:2016}. As it will be of importance in the current contribution, we mention that the curl of such a tangential continuous vector field exists in $L^2$ sense, which means that
\begin{align}
\opcurl_\Xv \uv \in L^2(\Omega) \label{curl}
\end{align}
can be used in a virtual work statement without constraints.

The final, linear variational formulation, as derived in \cite{PechsteinMeindlhumerHumer:2018}, reads: find $\uv \in \mathbf{V}(\Omega)$, $\taut \in \boldsymbol{\Sigma}(\Omega)$ and $\phi \in W(\Omega)$ such that for all admissible $\delta \uv$, $\delta \taut$ and $\delta \phi$ we have
\begin{align}
\begin{aligned}
- \int_\Omega (\Ctm^{-1}: \taut + \Ctm^{-1}:\Etm\cdot \nabla_\Xv \phi) : \delta \taut\, d\Xv + \langle \delta \taut, \epst(\uv)\rangle_\Omega + \langle \taut, \delta \epst(\uv)\rangle_\Omega & \\
- \int_\Omega (\Etm \cdot\Ctm^{-1}: \taut - (\Dtm - \Etm :\Ctm^{-1} :\Etm)\cdot \nabla_\Xv \phi)\cdot \delta \nabla_\Xv \phi\, d\Xv & = 0.
\end{aligned}
\end{align}

\subsection{A consistent linearization} \label{sec:mixedfe}

In the sequel, we derive the linearized update equation at some actual state $(\uv_0, \phi_0)$. We assume, that the displacement $\uv_0$ is continuous, i.e.\ there are no gaps between elements. In an implementation, this means that in each step we need to project the discontinuous displacement update to some continuous update. Moreover, we assume we are given the total stress $\taut_0$, which was computed in the last iterative step, in current configuraton. Note that $\taut_0$ may differ from $\overline{\taut}_0$, where $\overline{ \taut}_0$ is determined by $\uv_0$ and $\phi_0$ as the derivative of the free energy function $\Psi$. E.g., in the first order method, $\overline{\taut}_0$ is constant per element, while $\taut_0$ is linear per element and normal-normal continuous.

All finite element spaces are defined along \eqref{spaceX1} - \eqref{spaceX3}, but in the deformed mesh. This means we will find the displacement update $\uv_\Delta$, the symmetric part of the total stress $\taut^{sym}$ and the electric potential update $\phi_\Delta$, as well as the various virtual quantities, such that
\begin{align}
\uv_\Delta, \delta \uv &\in \mathbf{V}(\omega_0) := \{ \vv: \vv|_{t_0} \in P^{k}(t_0), \vv_t \text{ cont.}, \vv_t = 0\text{ on }\gamma_{0,fix}\},
\label{spacex1}\\
\taut^{sym}, \delta \taut & \in \boldsymbol{\Sigma}(\omega_0) := \{ \sigmat: \sigmat \text{ symm}, \sigmat|_{t_0} \in P^{k}({t_0}), \sigma_{nn} \text{ cont.}, \sigma_{nn} = 0\text{ on }\gamma_{0,free}\},
\label{spacex2}\\
\phi_\Delta, \delta \phi & \in W(\omega_0):= \{ w: w|_{t_0} \in P^{k+1}({t_0}), w \text{ cont.}, w = 0 \text{ on }\gamma_{0,E}\}.
\label{spacex3}
\end{align}
In finding a consistent linearized variational formulation of \eqref{pde1} - \eqref{pde4}, we have to overcome several problems. First, the update gradient $\Ft_\Delta$ does not exist as  such, but only in distributional sense. Second, the total stress $\taut $ is not symmetric in \eqref{pde1} due to the geometric stiffening effects, while the original TDNNS method uses the symmetry of the stress tensor as essential ingredient. We will discuss the implications and solution strategies in the following.

In a first step, we linearize \eqref{pde1} around $\uv_\Delta = 0,\taut = \taut_0$ and $\phi_\Delta = 0$, and obtain
\begin{equation}
\taut - \nabla_\xv \uv_\Delta \cdot \taut_0 = \overline{\taut}_0 + \ctm : \epst_{\Delta}  - \etm \cdot\nabla_\xv \phi_\Delta .\label{pde1lin}
\end{equation}
As mentioned above, in the mixed finite element method the displacement gradient $\nabla_\xv \uv_\Delta$ is available only as a distribution, and therefore needs special treatment. When we split the displacement gradient into a symmetric and a skew-symmetric part, we see that the symmetric part is the strain, while the skew-symmetric part corresponds to the curl of the displacement (compare \eqref{curlskw}), which is well-defined due to \eqref{curl}. We have in detail
\begin{align}
\nabla_\xv \uv_\Delta &= \epst_\Delta + \nabla^{skw}_\xv \uv_\Delta, \label{symskw}
\end{align}
Using equation \eqref{symskw} and basic linear algebra we can reformulate \eqref{pde1lin} as
\begin{align}
\taut - \overline{\taut}_0 - \nabla_\xv^{skw} \uv_\Delta \cdot \taut_0 + \etm^T \cdot \nabla_\xv \phi_\Delta &= \ctm : \epst_\Delta + \epst_\Delta \cdot \taut_0. \label{pde1lin2}
\end{align}

We split the nonsymmetric total stress $\taut$  as well as \eqref{pde1lin2} into a symmetric and a skew-symmetric part. Of course, at convergence the skew-symmetric part will tend to zero,
\begin{align}
\taut &= \taut^{sym} + \taut^{skw}, \label{tausymskw}\\
\taut^{sym} - \overline{\taut}_0 - \sym(\nabla^{skw}_\xv\uv_\Delta \cdot \taut_0) + \etm^T\cdot \nabla_\xv \phi_\Delta &= \ctm : \epst_\Delta + \sym(\epst_\Delta\cdot  \taut_0) =: \bar \ctm: \epst_\Delta, \label{pde1symlin2}\\
\taut^{skw} - \skw(\nabla^{skw}_\xv \uv_\Delta\cdot \taut_0) &= \skw(\epst_\Delta \cdot\taut_0). \label{pde1skwlin2}
\end{align}
In \eqref{pde1symlin2}, $\bar \ctm$ is implicitly defined as the stiffness tensor including (part of) the geometric stiffening effects due to $\taut_0$. Next,  we solve for $\epst_\Delta$ in \eqref{pde1symlin2},
\begin{align}
\epst_\Delta &= \bar\ctm^{-1}:(\taut^{sym} - \overline{\taut}_0 - \sym(\nabla^{skw}_\xv\uv_\Delta \cdot \taut_0)+ \etm ^T\cdot\nabla_\xv \phi_\Delta). \label{pde1symlin3}
\end{align}
We insert the result into \eqref{pde1skwlin2} to obtain an expression for $\taut^{skw}$,
\begin{align}
\taut^{skw} = \skw\Big(\!\big(\nabla^{skw}_\xv \uv_\Delta +  \bar \ctm^{-1}\!\!:(\taut^{sym} - \overline{\taut}_0 - \sym(\nabla^{skw}_\xv\uv_\Delta \!\cdot \!\taut_0)+ \etm^T\!\!\cdot \nabla_\xv \phi_\Delta)\big) \cdot\taut_0\Big). \label{pde1skwlin3}
\end{align}
Thus, we get the total stress $\taut$  by \eqref{tausymskw} as the sum
\begin{align}
\taut &= \taut^{sym}+ \skw\Big(\!\big(\nabla^{skw}_\xv \uv_\Delta +  \bar \ctm^{-1}\!\!:(\taut^{sym} - \overline{\taut}_0 - \sym(\nabla^{skw}_\xv\uv_\Delta \!\cdot \!\taut_0)+ \etm^T\!\!\cdot \nabla_\xv \phi_\Delta)\big) \cdot\taut_0\Big). \label{tau2}
\end{align}
Now we can gather the final set of partial differential equations. We use the material law \eqref{pde1symlin3} in the first line below. Then we insert the expression for $\taut$ from \eqref{tau2} in the mechanical balance equation and obtain \eqref{finalpde2}. Note that we used $\taut^{sym} = (\taut^{sym})^T$ and $\taut^{skw} = - (\taut^{skw})^T$. Additionally, we insert the expression for $\epst_\Delta$ from \eqref{pde1symlin3} in the electric material law \eqref{pde2} and further in Gauss' law \eqref{pde4} to gain \eqref{finalpde3}. From these equations below, we will derive the virtual work statement.
\begin{align}
- \bar \ctm^{-1}:(\taut^{sym} - \overline{\taut}_0 - \sym(\nabla^{skw}_\xv\uv_\Delta \cdot \taut_0)+ \etm^T\cdot \nabla_\xv \phi_\Delta)  + \epst_\Delta = 0,& \label{finalpde1}\\
\begin{aligned}
- \opdiv_\xv\!\big(\taut^{sym}\big)+\opdiv_\xv\!\big(\!\skw( \nabla^{skw}_\xv\uv_\Delta\cdot \taut_0) \big)  &\\
+ \opdiv_\xv\!\big(\!\skw((\bar \ctm^{-1}:(\taut^{sym} - \overline{\taut}_0 - \sym(\nabla^{skw}_\xv\uv_\Delta \cdot \taut_0)+ \etm^T\cdot \nabla_\xv \phi_\Delta))\cdot \taut_0)\big) &= 0, \end{aligned}\label{finalpde2}\\
- \opdiv_\xv\!\big(\overline{\dv}_0 - \etm :\bar\ctm^{-1}:(\taut^{sym} - \overline{\taut}_0 - \sym(\nabla^{skw}_\xv\uv_\Delta \cdot \taut_0)+ \etm ^T\cdot\nabla_\xv \phi_\Delta) + \dtm \cdot\nabla_\xv \phi_\Delta\big) = 0.& \label{finalpde3}
\end{align}
We multiply \eqref{finalpde1} with a symmetric virtual stress $\delta \taut$, \eqref{finalpde2} with a virtual displacement $\delta \uv$ and \eqref{finalpde3} with a virtual electric potential $\delta \phi$. We integrate over the deformed domain $\omega_0$ and add the equations, to arrive at
\begin{equation}
\begin{aligned}
-\int_{\omega_0} \delta \taut : \bar \ctm^{-1}:(\taut^{sym} - \overline{\taut}_0 - \sym(\nabla^{skw}_\xv\uv_\Delta \cdot \taut_0)+ \etm ^T\cdot\nabla_\xv \phi_\Delta)\, d\xv  +\langle \epst_\Delta, \delta \taut\rangle_{\omega_0}  & \\
- \langle \opdiv_\xv \taut^{sym}, \delta \uv\rangle_{\omega_0}  + \int_{\omega_0} \opdiv_\xv(\skw( \nabla^{skw}_\xv\uv_\Delta\cdot \taut_0) ) \cdot \delta \uv\, d\xv
& \\
+ \int_{\omega_0} \opdiv_\xv(\skw((\bar \ctm^{-1}:(\taut^{sym} - \overline{\taut}_0 - \sym(\nabla^{skw}_\xv\uv_\Delta\! \cdot\! \taut_0)+ \etm^T\!\cdot\! \nabla_\xv \phi_\Delta) )\cdot \taut_0) ) \cdot \delta \uv\, d\xv& \\
- \int_{\omega_0} \opdiv_\xv\!\big(\overline{\dv}_0 - \etm : \bar\ctm^{-1}\!:(\taut^{sym} - \overline{\taut}_0 - \sym(\nabla^{skw}_\xv\uv_\Delta \!\cdot\! \taut_0)+ \etm ^T\!\cdot\!\nabla_\xv \phi_\Delta)+ \dtm \cdot\nabla_\xv \phi_\Delta\big) \delta \phi\, d\xv & \\= 0.&
\end{aligned}
\end{equation}
Here we already used the distributional strain and divergence operators. Integration by parts and the equivalence \eqref{epsdiv1} - \eqref{epsdiv2}, together with boundary conditions, lead to
\begin{equation}
\begin{aligned}
-\int_{\omega_0} \delta \taut : \bar \ctm^{-1}:(\taut^{sym} - \overline{\taut}_0 - \sym(\nabla^{skw}_\xv\uv_\Delta \cdot \taut_0)+ \etm ^T\cdot\nabla_\xv \phi_\Delta)\, d\xv  +\langle \epst_\Delta, \delta \taut\rangle_{\omega_0}  & \\
+\langle \taut^{sym}, \delta \epst\rangle_{\omega_0}  - \int_{\omega_0}  (\nabla^{skw}_\xv\uv_\Delta\cdot \taut_0) : \delta \nabla^{skw}_\xv\uv\, d\xv
& \\
- \int_{\omega_0} ((\bar \ctm^{-1}:(\taut^{sym} - \overline{\taut}_0 - \sym(\nabla^{skw}_\xv\uv_\Delta \cdot \taut_0)+ \etm^T\cdot \nabla_\xv \phi_\Delta) )\cdot \taut_0) : \delta \nabla^{skw}_\xv\uv\, d\xv& \\
+ \int_{\omega_0} \big(\overline{\dv}_0 - \etm : \bar\ctm^{-1}\!:(\taut^{sym} - \overline{\taut}_0 - \sym(\nabla^{skw}_\xv\uv_\Delta \!\cdot\! \taut_0)+ \etm ^T\!\cdot\!\nabla_\xv \phi_\Delta)+ \dtm \cdot\nabla_\xv \phi_\Delta\big) \cdot\delta \nabla_\xv\phi\, d\xv &\\= 0.&
\end{aligned}
\end{equation}
We regroup terms and apply some more linear algebra for symmetric and skew-symmetric tensors. When collecting all affine linear terms on the right hand side, we obtain a symmetric, linear finite element formulation for $\uv_\Delta$, $\taut^{sym}$ and $\phi_\Delta$ in deformed configuration,
\begin{equation}
\begin{aligned}
-\int_{\omega_0} \delta \taut : \bar \ctm^{-1}:(\taut^{sym} + \etm ^T\cdot\nabla_\xv \phi_\Delta)\, d\xv  +\langle \epst_\Delta, \delta \taut\rangle_{\omega_0} +\langle \taut^{sym}, \delta \epst\rangle_{\omega_0} & \\
- \int_{\omega_0} \delta \nabla_\xv\phi \cdot ( \etm : \bar\ctm^{-1}:\taut^{sym}- (\dtm-\etm : \bar\ctm^{-1}:\etm^T) \cdot\nabla_\xv \phi_\Delta) \, d\xv &\\
+ \int_{\omega_0} (\delta \taut + \delta \nabla_\xv\phi \cdot \etm) : \bar \ctm^{-1}:\sym(\nabla^{skw}_\xv\uv_\Delta \cdot \taut_0)\, d\xv
& \\
+ \int_{\omega_0}\sym(\delta \nabla^{skw}_\xv\uv \cdot \taut_0) :\bar  \ctm^{-1}:(\taut^{sym} + \etm^T\cdot \nabla_\xv \phi_\Delta)  \, d\xv & \\
  - \int_{\omega_0}  (\nabla^{skw}_\xv\uv_\Delta\cdot \taut_0) : \delta \nabla^{skw}_\xv\uv\, d\xv- \int_{\omega_0}\sym(\delta \nabla^{skw}_\xv\uv \cdot \taut_0) : \bar \ctm^{-1}:\sym(\nabla^{skw}_\xv\uv_\Delta \cdot \taut_0)  \, d\xv &= \\
-\int_{\omega_0} (\delta \taut - \sym(\delta \nabla^{skw}_\xv\uv\cdot \taut_0)): \bar \ctm^{-1}: \overline{\taut}_0 \, d\xv    
- \int_{\omega_0}\delta \nabla_\xv\phi\cdot (\overline{\dv}_0 + \etm : \bar\ctm^{-1}: \overline{\taut}_0 ) \, d\xv . \\
\end{aligned}\label{eq:var}
\end{equation}

All terms in the above formulation are well-defined for TDNNS elements, as we need only the work pairs $\langle \taut^{sym}, \delta \epst\rangle_{\omega_0}$ and $\langle \delta \taut, \epst_\Delta\rangle_{\omega_0}$, the gradient of $\phi_\Delta$ and the curl of $\uv_\Delta$. The system matrix is then indefinite but symmetric.

\subsection{Nearly incompressible materials} \label{sec:incomp}

The method is applicable also for nearly incompressible materials with Poisson's ratio $\nu$ approaching $\tfrac{1}{2}$ or large Lam\'e parameter $\lambda$, see \cite[Chapter 5]{Sinwel:2009}. To ensure numerical stability, a consistent stabilization term is added to the variational form \eqref{eq:var},
\begin{align}
-\sum_{T \in \mathcal{T}}\int_T  h_T^2 \opdiv \taut^{sym} \cdot \opdiv \delta \taut\, d\xv.
\end{align}
Above, the parameter $h_T$ denotes the local mesh size, i.e.\ the diameter of the element. In case volume forces $\fv$ are present such that $-\opdiv \taut^{sym} = \frac{1}{J_0} \fv$, the right hand side of \eqref{eq:var} has to be modified accordingly by adding
\begin{align}
\sum_{T \in \mathcal{T}}\int_T \frac{h_T^2 }{J_0} \fv \cdot \opdiv \delta \taut\, d\xv.
\end{align}

\section{Numerical results} \label{sec:numerics}

\subsection{Compression of plate with circular hole}

\begin{figure}
\begin{center}
\includegraphics[width=0.9\textwidth]{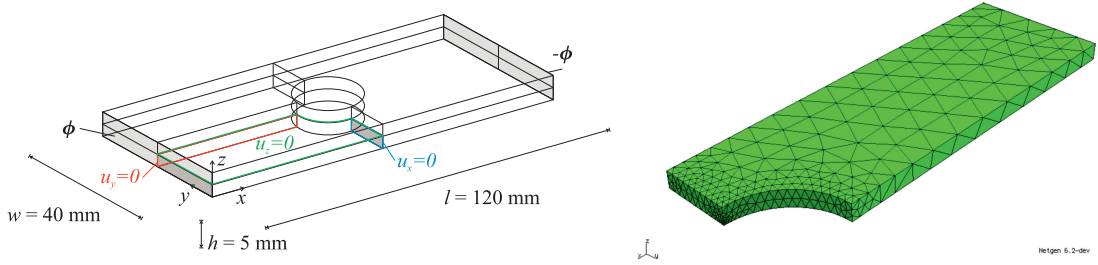}
\end{center}
\caption{Left: Problem definition of a plate with circular hole. Electrical boundary conditions and symmetry conditions are prescribed. Right: Tetrahedral mesh of one eighth of the geometry.} \label{fig:platehole}
\end{figure}

The first numerical example is taken from \cite{VuSteinmannPossart:2007}. We consider a plate of length 120 mm, width 40 mm and thickness 5 mm. The circular hole is centered and of radius 10 mm. The electrodes are located on both ends of the plate, such that the electric field points approximately in $x$ direction, see Figure~\ref{fig:platehole}. Due to symmetry, only one eighth of the geometry is meshed, and symmetric boundary conditions are imposed on the internal faces. 

We will assume that the electro-active polymer is described by a free energy function
\begin{align}
\Psi(\Ct, \Ev) &= \frac{\mu}{2} (\Ct : \It - d) - \mu \log J +  \frac{\lambda}{2} (\log J)^2 + c_1 \Ev\cdot \Ev + c_2 \Ev \cdot \Ct\cdot \Ev.
\end{align}
We neglect the vacuum permittivity $\epsilon_0$ due to its smallness. We use material parameters from \cite{VuSteinmannPossart:2007}, setting $\mu = 5$~N/mm$^2$, $\lambda = 20/3$ N/mm$^2$, $c_1 = 10$~N/V$^2$ and $c_2 = 6$~N/V$^2$. We apply the electric potential boundary conditions in load steps of $\Delta \phi = 10$~V each.

In each update step, the discontinuous displacement update $\uv_\Delta$ is interpolated by a continuous finite element function, which leads then to the new iterate $\uv_0$. The iteration is stopped as soon as the discontinuous finite element solution dominates the interpolated update by a factor of 10, as then the discretization error dominates the nonlinear residual. In the present example, we needed  5~Newton steps in each load increment to reach this convergence.

\begin{figure}
\begin{center}
\includegraphics[width=0.8\textwidth]{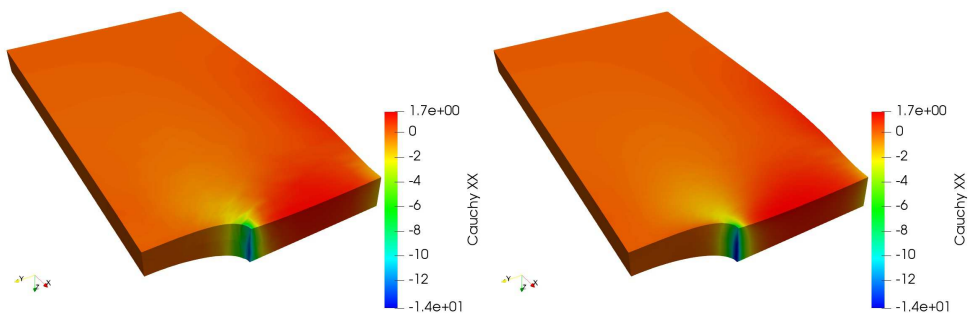}
\end{center}
\caption{Longitudinal total strain component $\tau_{xx}$ at $\Delta \phi = 2 \cdot 55$~V, left: polynomial order 1, right: polynomial order 2.} \label{fig:platehole_sxx}
\end{figure}
\begin{figure}
\begin{center}
\includegraphics[width=0.8\textwidth]{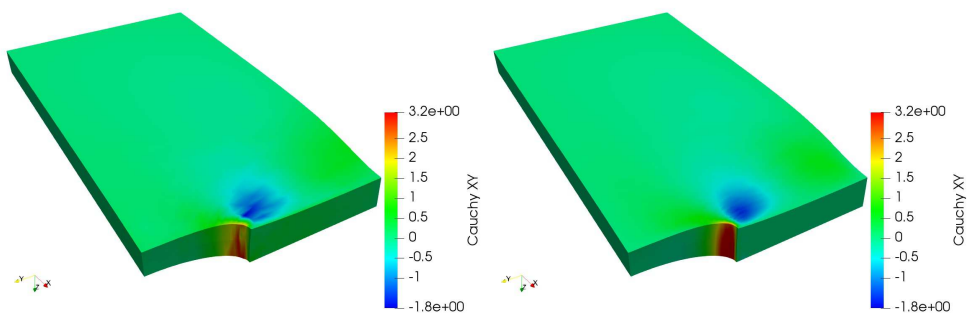}
\end{center}
\caption{Shear total strain component $\tau_{xy}$ at $\Delta \phi = 2 \cdot 55$~V, left: polynomial order 1, right: polynomial order 2.} \label{fig:platehole_sxy}
\end{figure}
\begin{figure}
\begin{center}
\includegraphics[width=0.8\textwidth]{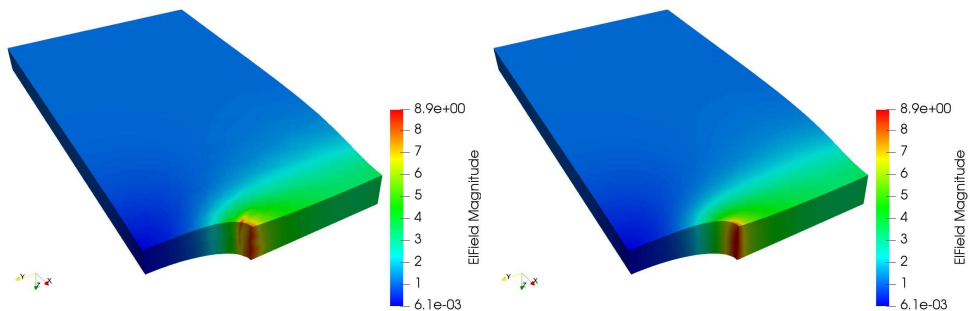}
\end{center}
\caption{Absolute value of electric field $|\Ev|$ at $\Delta \phi = 2 \cdot 55$~V, left: polynomial order 1, right: polynomial order 2.} \label{fig:platehole_Eabs}
\end{figure}

All computations are carried out in the framework of the open-source finite element package Netgen/NGSolve \cite{netgen}.
We use a tetrahedral mesh of  average mesh size 5~mm, while it is refined to 1~mm near the circular hole. The mesh consists of 1364 elements and is depicted in Figure~\ref{fig:platehole}. We compare results of the first and second order elements (i.e.\ the displacement and stress basis functions are of order one or two, respectively, while the electric potential basis functions are one order higher). This choice is implied by the convergence results from \cite{PechsteinMeindlhumerHumer:2018}. Electric field and total stress resulting for a potential difference of $\Delta \phi = 110$~V are depicted in Figures~\ref{fig:platehole_sxx} to \ref{fig:platehole_Eabs}. 

\subsection{Buckling of dielectric elastomer}

The second example is a  buckling actuator that was first presented in \cite{KlinkelZweckerMueller:2013} and also analyzed in \cite{StaudiglKrommerVetyukov:2017}. A thin membrane made of a dielectric elastomer is clamped laterally. Electrodes are applied to the top and bottom of the membrane. When a voltage is applied, the thickness of the  membrane decreases, while it extends in-plane. The restriction on the clamped edges leads to buckling behavior. The buckling is initiated by a small volumetric force acting in thickness direction. 

\begin{figure}
\begin{center}
\includegraphics[width=0.5\textwidth]{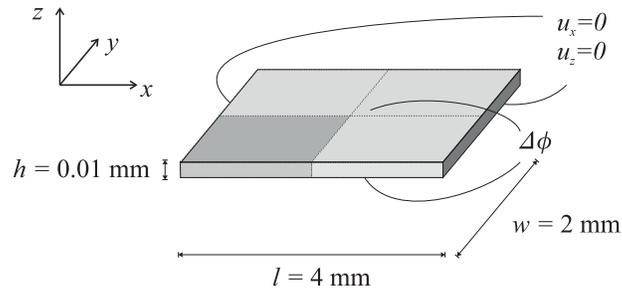}\
\end{center}
\caption{Problem definition for the buckling actuator. Due to symmetry, only the darker shaded area is discretized.} \label{fig:bucklingproblem}
\end{figure}

For the description of the dielectric elastomer we use a nearly incompressible Neo-Hookean material such as described in the latter reference \cite{StaudiglKrommerVetyukov:2017}. The free energy function of this dielectric elastomer is given by
\begin{align}
\Psi(\Ct, \Ev) = \frac{\mu}{2} (\operatorname{tr}(\Ct) - 3) - \mu \log J + \frac{\lambda}{2} (\log J)^2 - \frac12 \epsilon \Ev \cdot \Ct^{-1} \cdot \Ev.
\end{align}
with $\mu = 20\,689$ Pa, $\lambda = 100$ MPa and $\epsilon = (J + \chi)\epsilon_0$ with $\chi = 3.7$. We mesh a quarter of the geometry, applying symmetric boundary conditions at the interior boundaries. See Figure~\ref{fig:bucklingproblem} for a problem description.

\begin{figure}
\begin{center}
\includegraphics[width=0.8\textwidth]{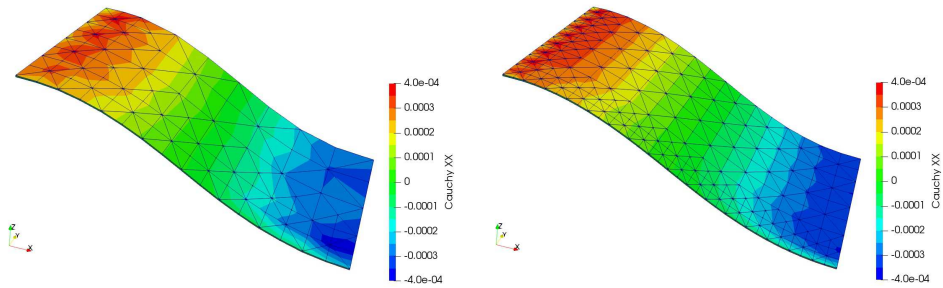}
\end{center}
\caption{Total stress $\tau_{xx}$ in the deformed (quarter) geometry for an applied voltage of $125$~V for the two different discretizations.} \label{fig:totalstress}
\end{figure}

\begin{figure}
\begin{center}
\includegraphics[width=0.8\textwidth]{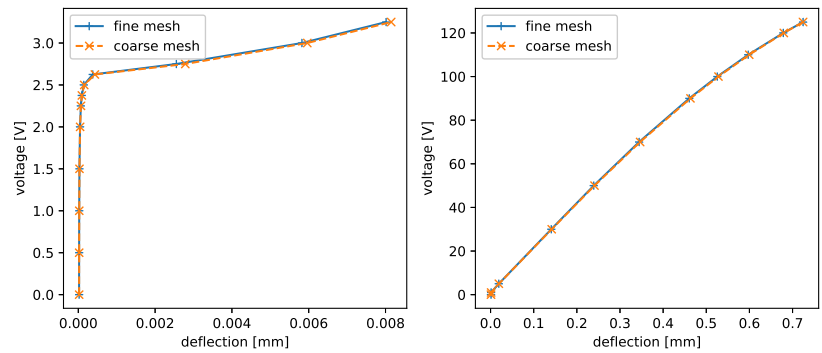}
\end{center}
\caption{Midpoint deflection versus applied voltage for the buckling actuator.} \label{fig:buckling}
\end{figure}

We use two different meshes, one consisting of 158 elements, the other one finer consisting of 458 elements. Both meshes are refined towards the (physical) boundaries of the quarter patch, to catch the singularities in the solution that arise there. In Figure~\ref{fig:totalstress} we plot the total stress $\tau_{xx}$ computed at $125$ V for both discretizations. Also in the buckling region, we could use quite a large update step from $1$ V to $5$ V, then we proceeded with even larger update steps as indicated in Figure~\ref{fig:buckling} to the right. Due to the near incompressibility of the material, we added the stabilization term from Section~\ref{sec:incomp}. Moreover, for stabilization, we modified the tangent compliance $\bar \ctm$ in the left hand side of \eqref{eq:var} by neglecting the geometric stiffening. Note that this change does not affect the solution, but only the convergence rate -- on average, 7 iterations were needed.

 In Figure~\ref{fig:buckling}, the midpoint deflection versus applied voltage is displayed for both meshes. The graphic on the left hand side provides the deflection  in the buckling area. The load is increased up to $125$~V, a corresponding evolution of the mid-point deflection can be found in the right hand side of Figure~\ref{fig:buckling}. The different mesh sizes lead to almost identical results. At $125$~V, the midpoint deflection amounts to 
$1.44866$~mm for the coares mesh and $1.44863$~mm for the fine mesh.
These values compare well to the value of $1.45$~mm reported in \cite{StaudiglKrommerVetyukov:2017}. In the original reference \cite{KlinkelZweckerMueller:2013}, where an Ogden-type material is used for the dielectric plate, the authors provide the value of $1.4$~mm for the opening.


\section{Conclusion}

In this work we have presented a mixed finite element formulation suitable for the description of large deformations in electro-active polymers. We assumed the material to be described by a free energy function, in our computational examples we used one of Neo-Hookean type. The element behaves well for nearly incompressible materials and is suitable for the discretization of slender structures by flat volume elements. Numerical examples are provided which support these claims. They show that Neo-Hookean materials with different electro-elastic coupling terms can be treated by the same formulation. The elements are freely available in the open-source software  Netgen/NGSolve {\tt ngsolve.org}.

The formulation is close to the updated Lagrangian algorithm. Therefore we hope that adaptive refinement or remeshing techniques, as were presented by \cite{LegerEtal:2014,LegerPepin:2016}, may be applicable in this framework. The possible application of these techniques, as well as the question of error estimation, shall be subject of further research.
Whether the elements are applicable to other formulations, which include e.g.\ multiplicative splitting of electric and elastic deformations as suggested in \cite{SkatullaSansourArockiarajan:2012,ZaehMiehe:2015}, may be of interest.

\bibliographystyle{plain}
\bibliography{Smart}

\end{document}